# LLAMA LIMA: A Living Meta-Analysis on the Effects of Generative AI on Learning Mathematics

Version 6, 11/26

Anselm Strohmaier, Samira Bödefeld, Oliver Straser, Frank Reinhold

University of Education Freiburg, Institute of Mathematics Education

**Abstract.**
The capabilities of generative AI in mathematics education are rapidly evolving, posing significant challenges for research to keep pace. Research syntheses remain scarce and risk being outdated by the time of publication. We present a Living Meta-Analysis (LIMA) on the effects of generative AI–based interventions on mathematics learning. We update the literature base every two months and apply a cumulative Bayesian multilevel meta-regression model. This paper reports results from the sixth version, including 37 studies, five of which were added since the fifth version. Results indicate a robust overall positive effect ($g = 0.52$, 95% CrI [0.32, 0.72]) and no evidence in favor of publication bias. Moderator analyses provided strong evidence for including integration mode, intervention dose, and grouping. Estimated effects were larger when generative AI supplemented teacher instruction, with greater intervention dose, and in group, teacher-led, or tutoring settings. For the remaining moderators, evidence was moderate or inconclusive.

**TLDR: Changes to the results.**
*More evidence:* The dataset increased to 37 studies, 74 effect sizes, and 6,061 participants. Five new studies contribute nine effects; two previously included studies were removed.
*Consolidated effect:* Hedges' g decreased from 0.57 to 0.52.
*Additional moderator findings:* In addition to previously reported strong evidence for a stronger effect when generative AI supplemented teacher instruction, this version reports strong evidence for larger effects at greater intervention doses, and in group, teacher-led, or tutoring settings than in individual settings. Moderate evidence favored including educational level and output mode; estimates were higher in primary education and with picture or video output.

**Living Meta-Analysis Notice**

LLAMA LIMA is maintained as a versioned living meta-analysis, with each update constituting a new version. This document is Version 6, 11/26. Updated versions of the meta-analysis will be published every two months and will be made publicly available through arXiv until LLAMA LIMA is retired and published permanently. Readers are advised to consult the most recent version when citing or interpreting the findings. This publication has not yet undergone formal peer review. The authors welcome feedback, critique, and suggestions.

**Living status:** Ongoing. The date of the last literature search for this version was October 1, 2026. The next literature search is scheduled for December 2026. The next version is scheduled for January 2027.

**Methodological changes since the previous version:** Smaller prior standard deviations are used for categorical moderators. Intervention dose (in hours, log-transformed for analyses) is now used as a moderator instead of intervention span. Grouping, output mode, research design, and assessment timepoint were added as new moderators. Whether prior AI literacy and general abilities were assessed was coded but not analyzed as a moderator. Publication type now distinguishes Q1 journal articles (by JCR) from other journal articles.

## 1 Introduction

Recent advances in *large language models* (LLMs), such as ChatGPT, have led to the rapid uptake of *generative artificial intelligence* (generative AI) in educational practice and research (Kasneci et al., 2023; Ng et al., 2025). Generative AI refers to AI systems, including LLMs, that generate novel content based on patterns learned from large-scale training data.

AI-based technologies have long been investigated in education (e.g., intelligent tutoring systems, adaptive learning, automated feedback; Chiu et al., 2023; Holmes et al., 2019), but generative AI has accelerated and shifted these efforts. Compared with many earlier technologies, generative AI enables open-ended natural-language interaction, allowing for new forms of instructional support such as spontaneously generating explanations, providing adaptive feedback, supporting problem solving, or facilitating individual and collaborative learning; at the same time, these systems typically offer limited predictability of outputs and limited transparency about how outputs are generated, giving rise to concerns regarding accuracy, interpretability, bias, and instructional control (Kasneci et al., 2023).

The rapid uptake of generative AI has been accompanied by a rapidly growing body of experimental research and research syntheses. Meta-analyses are particularly influential because they condense complex and heterogeneous evidence into apparently straightforward estimates of effectiveness that can inform research, practice, and policy. In this setting, the credibility of research syntheses is particularly important, and a growing body of research documents quality problems in this rapidly expanding literature. In a recent forensic audit, O'Neill (2026) examined 14 peer-reviewed meta-analyses on AI in education, some of which received considerable attention and are widely cited, and found serious quality issues in all of them. Fütterer et al. (2026) similarly identified limitations in the quality, reporting, and reproducibility of early syntheses and their primary evidence.

O'Neill (2026) criticized research that asks whether “AI” or “ChatGPT” improves learning, treating a technology as though it were a coherent intervention. Educational effects instead depend on the designed learning activity in which generative AI is used, the role assigned to the system, the learners and content involved, and the outcome being measured. This also implies that the benefits of generative AI are unlikely to be uniform across educational domains. Mathematics learning is characterized by domain-specific cognitive processes and instructional practices (Schneider & Stern, 2010; Schoenfeld, 2016). Consequently, both the potential benefits and the challenges associated with implementing generative AI in mathematics education may differ from those observed in other domains (Cho, 2026; Pepin et al., 2025; Walkington, 2025).

Within mathematics education, consolidated and accountable empirical evidence that accounts for the specific characteristics of the field is still lacking (Cho, 2026; Fock & Siller, 2025a, 2025b; Pepin et al., 2025; Su et al., 2026; Turmuzi et al., 2026; Walkington, 2025). In particular, Su et al. (2026) point out that studies on standardized mathematical learning outcomes remain scarce, which limits the conclusions that can be drawn about learning effects. Accordingly, two central questions remain open: (a) whether generative AI–based instructional interventions can specifically support mathematics learning effectively, and (b) which factors determine their effectiveness.

While both questions can be investigated with meta-analyses, the rapid evolution of generative AI technologies—and, with them, the research landscape—means that conventional research syntheses risk being unable to keep up, and serious concerns exist about the credibility of hastily curated syntheses. We therefore address these questions through LLAMA LIMA (*Large Language Models and Generative AI in Mathematics Education: A Living Meta-Analysis*), a living meta-analysis in which we continuously expand our database, update our analyses, and make the results immediately accessible. Moreover, we openly invite feedback and critique and transparently document changes and their reasons across versions. LLAMA LIMA was first published in January 2026, and it has been updated every two months since then. To our knowledge, this represents the first publication-based living meta-analysis in educational research.

### 1.1 The Challenge of Synthesizing Experimental Research on Generative AI in Mathematics Education

Although generative AI can potentially support mathematics learning in various ways (Fock & Siller, 2025a, 2025b; Pepin et al., 2025; Su et al., 2026; Turmuzi et al., 2026; Walkington, 2025), the existing empirical research is highly heterogeneous, encompassing different instructional designs, pedagogical purposes of AI, outcome measures, and learner populations. Many studies rely on exploratory designs, short-term interventions, and narrowly defined instructional contexts, which limits the generalizability of individual findings. Under these conditions, individual studies are often insufficient to provide robust answers to questions about overall effectiveness, despite the growing use of generative AI in educational practice and the strong interest of researchers, educators, policymakers—and learners. Systematic research syntheses and meta-analyses help integrate and summarize findings across studies.

However, the rapid technological development of generative AI poses a challenge to the synthesis of existing evidence. Conducting and publishing rigorous empirical studies—and summarizing them in conventional meta-analyses—requires substantial time (Huisman & Smits, 2017). This is illustrated by the fact that many existing reviews on AI in education predate the widespread availability of large language models (e.g., Tlili et al., 2025), which began with the release of ChatGPT in late 2022. More recent syntheses illustrate how quickly evidence assessments become outdated as model capabilities evolve: For example, the scoping review by Pepin et al. (2025), published in February 2025 and based on studies available up to May 2024, discusses limitations in ChatGPT's mathematical performance that have mostly been mitigated by subsequent model versions. Similarly, the meta-analysis by Fock and Siller (2025b), published in November 2025 and based on a literature search conducted in November 2024, identified only three experimental studies in secondary mathematics education explicitly investigating generative AI. More recently, the meta-analysis by Liu et al. (2026), which includes studies published up to September 2025, identified 22 eligible studies. However, it also incorporated theses and other non–peer-reviewed literature and, in some cases—contrary to its stated inclusion criteria—studies without control groups. The recent meta-analysis by Boolzen et al. (2026) on effects of generative AI in STEM domains includes only four studies in mathematics up to May 7, 2026, but notes that the number of studies is growing quickly.

Another challenge concerns trustworthiness and accountability. O'Neill (2026) shows that this concern extends beyond isolated errors: broad meta-analyses can combine incoherent constructs, weakly vetted primary studies, and inadequate statistical safeguards while still becoming highly influential. One of the most highly cited meta-analyses on the effects of ChatGPT on learning was the paper by Wang and Fan (2025), in which the authors reported substantial positive effects. Only six days later, Bartoš, Martinková, et al. (2025) published a reanalysis of the same data indicating substantial publication bias, receiving only a fraction of the attention. Almost one year later, on April 22, 2026, the editors retracted the Wang and Fan (2025) paper. The retraction note states that there were “discrepancies in the meta-analysis” (p. 1), undermining “the confidence the Editor can place in the validity of the analysis and the resulting conclusions” (p.1, Wang & Fan, 2026). Arguably, in a rapidly evolving field like that of generative AI in education, mistakes will be made at the level of individual studies as well as in research syntheses. However, this becomes highly problematic when the published record remains frozen at the point of initial publication, and no established practice or editorial process exists for correcting or updating results or excluding retracted studies (O'Neill, 2026).

These examples illustrate that conventional research syntheses may no longer fully reflect the state of the technology when they become available and offer limited mechanisms for scrutiny, corrections, and updates. This mismatch between rapid technological development and the availability of timely, trustworthy research syntheses poses a substantial methodological challenge.

### 1.2 Living Meta-Analyses for Rapidly Evolving Evidence

While this sense of urgency and the volatility of the research landscape are relatively uncommon in educational sciences, other domains are more familiar with rapidly evolving evidence. In response, the field of medicine has developed approaches known as *living evidence syntheses* (Elliott et al., 2014): *Living systematic reviews* (LSRs) and *living meta-analyses* (LIMAs) are designed to reduce the time lag of conventional research syntheses by enabling continuous updating of the evidence base.

Living evidence syntheses differ from conventional research syntheses in several key respects:

- Continuous publication: Rather than being published as one-time products, living evidence syntheses are designed to be updated regularly, allowing the synthesis to reflect the evolving state of the evidence and to include corrections and adaptations.
- Repeated literature searches: Searches are conducted at regular intervals, often supported by automated procedures, and newly published studies are screened, coded, and incorporated with minimal delay.
- Inclusion of recent evidence: To reduce publication lag, living evidence syntheses may include preprints and conference proceedings alongside peer-reviewed studies. This approach involves a trade-off between rapid inclusion of new evidence and traditional publication quality filters.

- Bayesian statistics: Bayesian meta-analytic models allow cumulative statistical modeling, with estimates updated as the evidence base expands. Moreover, Bayesian model comparisons allow insufficient evidence for an effect to be distinguished from evidence for the absence of an effect.
- Correction and accountability: Living evidence syntheses can respond when primary studies are corrected or retracted, authors provide additional data, or coding and analytic decisions require revision. Version histories make these changes and their consequences transparent.

Together, these features make living meta-analyses particularly suitable for rapidly evolving research areas such as generative AI in education. A living format does not guarantee validity, but it allows the evidence base to remain current, transparent, and correctable. Traceable study and effect-size records also enable updates following corrections or retractions.

Meaningful synthesis also requires a clear account of the interventions being synthesized. O'Neill (2026) criticizes meta-analyses that treat AI tools as coherent interventions while pooling across substantial differences in instructional purpose, implementation, learner populations, and outcomes. For mathematics education, we therefore first distinguish the purposes of generative AI use and then consider further characteristics that may explain differences in learning effects.

### 1.3 Purposes of Generative AI Use in Mathematics Education

To structure the purposes for which generative AI may support mathematics learning, we draw on recent literature reviews and frameworks of generative AI integration in mathematics education (Chao & Han, 2026; Pepin et al., 2025; Schindler et al., 2026; Su et al., 2026; Turmuzi et al., 2026; Walkington, 2025). We propose a set of five categories that describe potential purposes for the use of generative AI to support students' mathematical learning.

**Direct mathematical support for learners.** Generative AI systems can generate correct answers and complete solutions for a wide range of school-relevant mathematical tasks (e.g., Getenet, 2024; Plevris et al., 2023; Strohmaier et al., 2025; Udias et al., 2024). This expert-like capability may support mathematics learning by allowing students to check their solutions and compare alternative solutions, much as they would with a calculator offering extended capabilities. This category covers learners obtaining mathematical answers, explanations, checks, or feedback but without a particular dialogic mode.

**Learner-AI dialogue.** Generative AI systems, and chatbots in particular, can be used by students as dialogic partners in multi-turn conversations. This can range from unguided, open conversations to structured, quasi-tutoring interactions. During extended interactions, generative AI can consider students' questions, solutions, and learning processes and provide individualized and adaptive assessment, feedback, and learning pathways (e.g., Chiu et al., 2023; Giannakos et al., 2025; Kasneci et al., 2023; Pepin et al., 2025; Wardat et al., 2023).

**Generation of learning materials**. Generative AI can be used by students to generate reusable learning materials, such as written or animated explanations, tasks, and exercises (e.g., Pepin et al., 2025; Schorcht et al., 2024; Turmuzi et al., 2026; Wu et al., 2024). When informed by typical learning trajectories, conceptual ideas and misconceptions, as well as information about the mathematical content and associated educational theory, generative AI can provide materials that are closely aligned with specific groups of learners, curricular goals, grading criteria, and assessment requirements (Turmuzi et al., 2026).

**Facilitation of collaborative learning.** Generative AI may further support mathematics learning by facilitating collaborative activities, such as group problem solving or project-based work (Pepin et al., 2025; Turmuzi et al., 2026). In these contexts, AI systems can be used to help structure collaboration, support discussion around mathematical ideas, or scaffold joint work on tasks.

**Teacher- or researcher-mediated support.** Finally, generative AI may support mathematics learning indirectly by assisting teachers or researchers. This includes preparing instructional content for learners, providing mathematical assistance, engaging in dialogue, or creating materials for professional users. Learners receive the resulting instruction or materials without interacting with the AI themselves. Buchholtz and Huget (2024) illustrate that systems such as ChatGPT can be used to generate lesson plans and provide feedback to teachers. To the extent that such uses contribute to higher-quality instruction, they may translate into improved student learning outcomes (Giannakos et al., 2025).

### 1.4 Determinants of the Effectiveness of Generative AI in Mathematics Learning

The purpose of AI use represents one dimension along which interventions using generative AI in mathematics differ. For many other study characteristics, theoretical considerations imply that effects on learning mathematics might differ. The following overview introduces potential moderators, grouped into AI-specific substantive characteristics, general substantive characteristics, and methodological and extrinsic characteristics (Lipsey & Wilson, 2000).

**AI-specific substantive characteristics.** Interventions using generative AI differ substantially in their design and implementation. Relevant characteristics that are specific to the implementation of generative AI include the social role of the AI, the intended purpose (see Section 1.3), the degree of learner autonomy, and whether generative AI is used as a supplement to, or replacement for, existing instruction, which we refer to as the *integration mode*. In addition, both the technological capabilities of the AI system and its output mode (e.g., text, pictures, or videos) could influence learning. Such design features are central to theories of technology-enhanced learning (e.g., Reinhold et al., 2024) and are likely to account for variability in observed effects across studies.

**General substantive characteristics.** Effects of generative AI–based interventions may vary as a function of general characteristics of the intervention that are not specific to learning with generative AI. Research on educational technologies highlights the importance of aligning instructional support with learners' prior knowledge and abilities (Hillmayr et al., 2020). Accordingly, learner characteristics are considered as potential moderators to examine whether effects differ across student populations. Beyond individual learners, the educational context in which generative AI is used may moderate its effectiveness. Relevant contextual factors include the mathematical content addressed: For example, understanding and addressing word problems may be a more natural task for LLMs, whereas geometry problems require different technological capabilities and non-linguistic mathematical concepts, such as handling visual representations, which could make both solving the problem itself and transferring it into an educational context a challenging task for generative AI. Moreover, the grouping of students during the intervention (e.g., individual, collaborative, or teacher-guided learning) could influence how effectively generative AI can support learning. Finally, total exposure time may influence the effectiveness of an intervention (Hillmayr et al., 2020). Another characteristic that might moderate the effectiveness of generative AI interventions is the underlying theory of learning guiding their design. Across studies, generative AI may be embedded within different instructional paradigms—such as direct instruction, problem-based or inquiry learning, cooperative or peer learning, metacognitive scaffolding, or conceptual change approaches—which shape how learners interact with the system and what kinds of cognitive processes are supported (e.g., Hattie & Yates, 2014).

**Extrinsic and methodological characteristics.** Finally, the effects of an intervention can be related to extrinsic and methodological characteristics. These characteristics are usually of less practical educational relevance but help situate and interpret the findings. They include publication type (e.g., proceedings, book chapters, or journal articles), research setting (e.g., laboratory vs. field), research design (e.g., individually randomized, cluster-randomized, or non-randomized), and post-test assessment timepoint (e.g., immediate vs. delayed).

## 2 The Present Study

Building on the potential of generative AI for mathematics learning, the fragmentation of empirical findings, and the limitations of conventional, static research syntheses in rapidly evolving fields, the present study adopts a living meta-analytic approach to integrate the available evidence. LLAMA LIMA (*Large Language Models and Generative AI in Mathematics Education: A Living Meta-Analysis*) aims to provide an up-to-date and continuously expanding synthesis of intervention studies examining the effects of generative AI on mathematics learning outcomes. In line with suggestions for living research syntheses outlined by Elliott et al. (2014), we adopt frequent updates, regular literature searches, the inclusion of preprints, and a Bayesian multilevel meta-analytic approach. We follow established reporting guidelines (PRISMA-LSR; Akl et al., 2024). LLAMA LIMA estimates a pooled effect across eligible interventions as an overall description of the accumulated evidence. Because interventions using generative AI vary in their design and implementation, the main interpretive focus lies on moderator analyses examining how effects differ across these conditions.

Accordingly, LLAMA LIMA addresses the following research questions:

**RQ1:** What is the pooled effect of instructional interventions using generative AI on mathematics learning outcomes?

**RQ2:** Which study characteristics influence the effectiveness of generative AI–based interventions for learning mathematics?

## 3 Methods

### 3.1 Living Mode Parameters

The planned schedule of this living meta-analysis is to update the literature search every two months and publish the revised analysis in the following month (i.e., the next database search is scheduled for December 2026, and the next version is scheduled for January 2027). Depending on the frequency of new publications and their influence on the overall effect and the feasibility of moderator analyses, these intervals might be altered in the future.

Reports excluded in previous versions may be included in subsequent versions if (a) previously unavailable reports become accessible or (b) authors contacted because of insufficient reported data provide additional information.

We plan to eventually end the living mode and publish a final version, but no timeline has been prespecified.

### 3.2 Literature Search

The literature search for the current version of the living meta-analysis was conducted on October 1, 2026, using Scopus for documents and preprints. The search targeted experimental and quasi-experimental studies on the use of generative AI in mathematics education and used a predefined combination of terms related to generative AI, education, mathematics, and intervention designs. In addition, a semantic search was conducted on OpenAlex. The full search strings are available in Appendix A.1.

Following title and abstract screening, full-text eligibility assessment, and removal of duplicate reports, five new studies met the inclusion criteria and were included in the current version. Details of the study identification, screening, eligibility assessment, and inclusion process across all versions are summarized in the PRISMA-LSR flow diagram shown in Figure 1.

During screening, studies were included if they (a) reported original data from an experimental or quasi-experimental intervention study, (b) included generative AI interaction in the intervention and no generative AI interaction in the control group, (c) involved human learners, (d) reported mathematics performance as an outcome measure, and (e) were written in English. We included studies published in peer-reviewed journals, edited book chapters, and conference proceedings, as well as preprints that were undergoing peer review at the time of screening. Theses, dissertations, reports, and other unpublished materials not subject to peer review were not included.

During full-text eligibility assessment, we excluded reports that (a) duplicated studies already represented in the database, retaining the most recent eligible report, (b) did not report sufficient data and whose authors had not responded to our request for additional data by the time of publication, or (c) did not have a research design suitable for meta-analysis. Studies were also excluded at this stage if the full text revealed that they did not meet the initial inclusion criteria. Finally, studies that were internally inconsistent or potentially not trustworthy (e.g., because of non-existent references or statistically implausible or contradictory data) were only included after a two-step review process. First, all authors independently and thoroughly evaluated the study. If concerns remained, the authors of the study were contacted and asked for clarification. If no response was received, the study was excluded. Exclusions due to missing information, inconsistencies, or trustworthiness are preliminary. If we receive additional information, these studies will be re-evaluated for future versions.

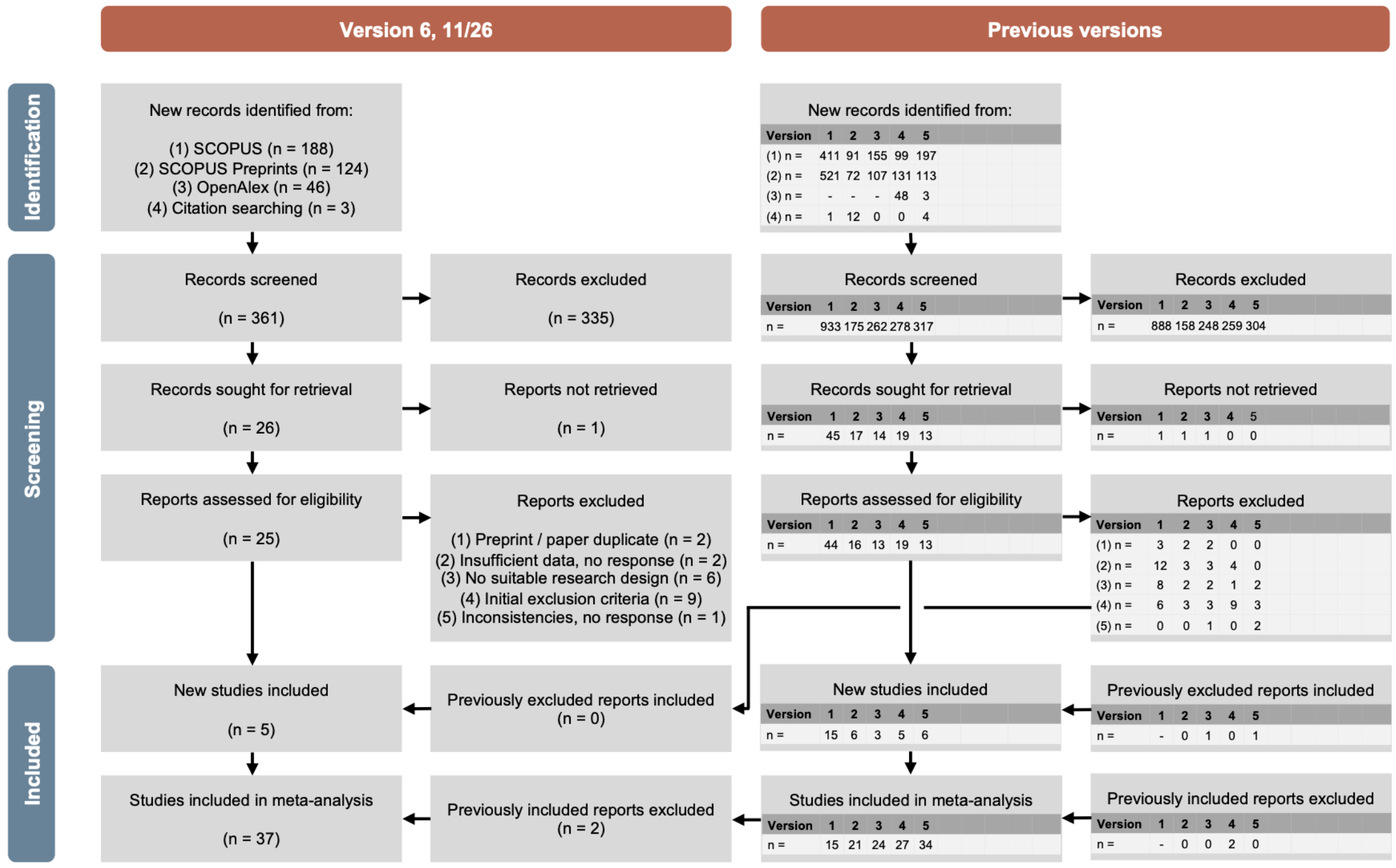

**Figure 1.** PRISMA-LSR Flow diagram documenting the study selection process for the living meta-analysis.

### 3.3 Moderator Coding

To analyze which factors influence whether generative AI affects mathematics learning, we continuously develop a moderator framework and coding manual following the guidelines by Lipsey and Wilson (2000). Moderators are grouped into general intervention, AI-specific intervention, methodological, and extrinsic characteristics. All studies are coded by two authors using a predefined manual; discrepancies are resolved by consensus. Thirteen moderators are reported in the current version.

**General intervention characteristics.** These included learners' educational level based on the *International Standard Classification of Education* (ISCED; Unesco, 2012), and information on prior mathematical knowledge, general abilities, and AI familiarity. Of these learner characteristics, only educational level was analyzed; prior general abilities were rarely assessed, and prior AI literacy was not assessed in any of the included studies. Mathematical content was categorized according to the NCTM standards (Number and Operations, Algebra, Geometry, Measurement, Data Analysis and Probability; NCTM, 2000). Grouping distinguished individual, group, teacher-led, and tutoring settings. Both content and grouping were included as moderators. Intervention dose was coded as total session time in hours, capturing exposure more directly than elapsed calendar span. Dose was analyzed continuously after log transformation and standardization. Estimates are reported at one hour and per doubling of dose.

**AI-specific intervention characteristics.** These included the purpose of AI use (see Section 1.3), modification of the tool or its output to improve mathematical or pedagogical quality, output mode (text, pictures, video), and state of technology, based on the reported model or version or, when this information was unavailable, the intervention year. For integration mode, we distinguished four categories: S*tandalone* means that no teacher was involved in either group. *Replacement* means that the control group was instructed by a teacher, whereas the intervention group was not. *Supplement* means that a teacher was involved in both groups. S*upported* means that a teacher was involved only in the intervention group.

**Methodological characteristics.** These included research design (randomized, cluster-randomized, non-randomized), research setting (field, controlled environment, laboratory), and assessment timepoint (immediate, delayed). Sample sizes and control conditions were also recorded.

**Extrinsic characteristics.** Publication type distinguished Q1 journal articles (by JCR), other journal articles, edited book chapters, conference proceedings, and preprints. Journal quartile was treated as a publication characteristic rather than an assessment of primary-study quality. Publication date and journal focus were also recorded.

### 3.4 Effect Size Extraction

For each study, we calculated standardized mean differences (Hedges' *g*) in outcomes between the intervention and control groups, using small-sample corrections. If there were several control groups (i.e., groups that did not interact with generative AI), we selected the control group with the most similar instructional activities and disregarded other control groups. When pre-test scores were available and on the same scale as the post-test, we calculated *g* as the between-group difference in gain scores, standardized by the pre-test standard deviation (Borenstein et al., 2009). If pre-test scores were not available but group equivalence was reported or implied by design (e.g., propensity score matching), post-test mean differences, standardized by the post-test standard deviation, were used. When means or standard deviations were unavailable, we transformed effect sizes from reported statistics. All effect sizes were coded such that positive values indicate a benefit of the intervention using generative AI.

Sampling standard errors for each effect were estimated according to the procedures described by Borenstein et al. (2009). To appropriately model the covariance of sampling errors within studies, a full variance–covariance matrix of the sampling errors was estimated in accordance with Viechtbauer (2010) assuming a correlation between measurements within the same group of $\rho = .7$ and an autocorrelation for several measurements across timepoints of $\varphi = .8$ (sensitivity analyses with $0 \leq \rho \leq .9$ and $0 \leq \varphi \leq .9$ revealed no practical influence of the choice of these parameters on the overall effect, $0.516 \leq g \leq 0.536$).

## 3.5 Data Analysis

**Bayesian multilevel meta-analysis.** To estimate the overall effect of generative AI–based instructional interventions on mathematics learning outcomes, we fitted a Bayesian multilevel meta-analytic model (Harrer et al., 2021). Main analyses were conducted in R using the *brms* package (Bürkner, 2017). The model accounted for the hierarchical structure of the data by treating effect sizes as nested within studies. In a Bayesian analysis, each parameter is represented by a posterior probability distribution; reported point estimates correspond to the posterior median, and uncertainty is summarized by 95% credible intervals. Residual heterogeneity was summarized separately at the study and effect-size levels using the posterior medians and 95% credible intervals of the corresponding random-effect standard deviations ($\tau$).

**Moderator analysis.** To quantify the amount of heterogeneity accounted for by a moderator, we calculated a median-based proportional reduction in heterogeneity, $R^2$, by comparing the posterior median of the heterogeneity variances from the baseline model without the moderator, $m_0$, with the corresponding posterior median from the moderator model, $m_1$. As an index of evidence for including each moderator, we computed Bayes factors ($BF_{10}$) based on marginal likelihoods estimated via bridge sampling (Gronau et al., 2020). Bayes factors greater than 1 quantify the relative evidence in favor of the moderator model, whereas values below 1 indicate evidence in favor of the baseline model. In line with common categories (Lee & Wagenmakers, 2014), we interpret Bayes factors between 1 and 3 (and between 0.33 and 1) as *anecdotal* evidence, between 3 and 10 (and between 0.10 and 0.33) as *moderate* evidence, between 10 and 100 (or between 0.01 and 0.10) as *strong* evidence, and values greater than 100 (or smaller than 0.01) as *extreme* evidence. Moderators were fitted separately to limit overfitting given the number of independent studies and sparse categories. These analyses therefore describe associations without disentangling correlated study characteristics. As evidence accumulates, multivariate meta-regression analyses will be considered.

**Prior specification and sensitivity analyses.** Weakly informative priors were specified for all model parameters. The pooled average effect was assigned a normal prior, $\mu \sim N(0, 1^2)$, and between-study and within-study heterogeneity standard deviations were assigned exponential priors, $\tau \sim \mathrm{Exp}(1/0.30)$. Categorical moderators were coded using normalized Helmert contrasts, with coefficient priors $\beta_j \sim N\big(0, 0.04 \cdot K/(K-1)\big)$, where $K$ is the number of categories. This scaling keeps the prior standard deviation of category effects at 0.20 across moderators. The coefficient for the standardized continuous moderator was assigned $\beta_j \sim N(0, 0.3^2)$. Because the predefined variance–covariance matrix of sampling errors was incorporated directly into the model, the residual variance parameter was fixed at one.

Prior sensitivity was examined for both the overall-effect and moderator analyses. For the pooled effect, heterogeneity, and categorical moderators, three prior specifications were compared: the default specification, a more informative specification, and a less informative specification. Relative to the default priors, the more informative specification halved the prior scales, using $\mu \sim N(0, 0.5^2)$, $\beta_j \sim N\big(0, 0.01 \cdot K/(K-1)\big)$, and $\tau \sim \mathrm{Exp}(1/0.15)$. The less informative specification increased the corresponding scales, using $\mu \sim N(0, 2^2)$, $\beta_j \sim N\big(0, 0.09 \cdot K/(K-1)\big)$, and $\tau \sim \text{Exp}(1/0.6)$. Thus, smaller exponential means represented more informative priors that concentrated more strongly near zero, whereas larger means represented less informative priors that assumed greater heterogeneity. Two additional sensitivity analyses used positively centered priors for the pooled average effect, $\mu \sim N(0.2, 1^2)$ and $\mu \sim N(0.4, 1^2)$, while retaining the default moderator and heterogeneity priors. Robustness was evaluated based on variation in the posterior median pooled effect (for overall-effect estimation) and Bayes factors (for moderator analyses) across prior specifications. For the continuous moderator, sensitivity analyses compared normal priors with mean zero and standard deviations of 0.15, 0.30, and 0.60.

Across prior specifications, the posterior median pooled effect ranged from $g$ = 0.50 to $g$ = 0.52. Among the reported moderators, the direction of evidence remained unchanged across prior specifications for all but research design, for which Bayes factors remained close to 1. For grouping and output mode, however, the strength of evidence was not stable. Bayes factors consistently favored moderator inclusion but ranged from $BF_{10}$ = 3.22 to 84.20 for grouping and from $BF_{10}$ = 1.53 to 7.48 for output mode. Thus, interpretations ranged from moderate to strong evidence for grouping and from anecdotal to moderate evidence for output mode, depending on the prior specification.

**Model estimation and inference.** Model parameters were estimated using Markov chain Monte Carlo sampling. Four independent chains were run with 2,000 warm-up iterations and 10,000 post-warm-up sampling iterations per chain for main analyses and prior sensitivity analyses, and 1,000 warm-up iterations and 3,000 post-warm-up sampling iterations per chain for ρ and φ sensitivity analyses and plots. Convergence was assessed using $\hat{R}$ and effective sample sizes, which indicated satisfactory convergence. Model fit was evaluated using visual checks of density overlay plots, which indicated adequate correspondence between observed and model-implied distributions. Results are reported as posterior medians accompanied by 95% credible intervals. These intervals represent the range within which the parameter lies with 95% posterior probability, given the observed data and model assumptions.

# 4 Results

## 4.1 Descriptive results

The meta-analysis included 37 studies contributing 74 effect sizes. Overall, these studies included data from 6,061 participants. Table A.2 in the appendix shows the included studies and basic study characteristics.

## 4.2 Pooled effect size

The posterior median pooled effect was $g$ = 0.52, with a 95% credible interval of [0.32, 0.72], indicating a positive effect accompanied by substantial uncertainty. Sensitivity analyses using different priors did not result in substantial variation in the estimates. Between-study

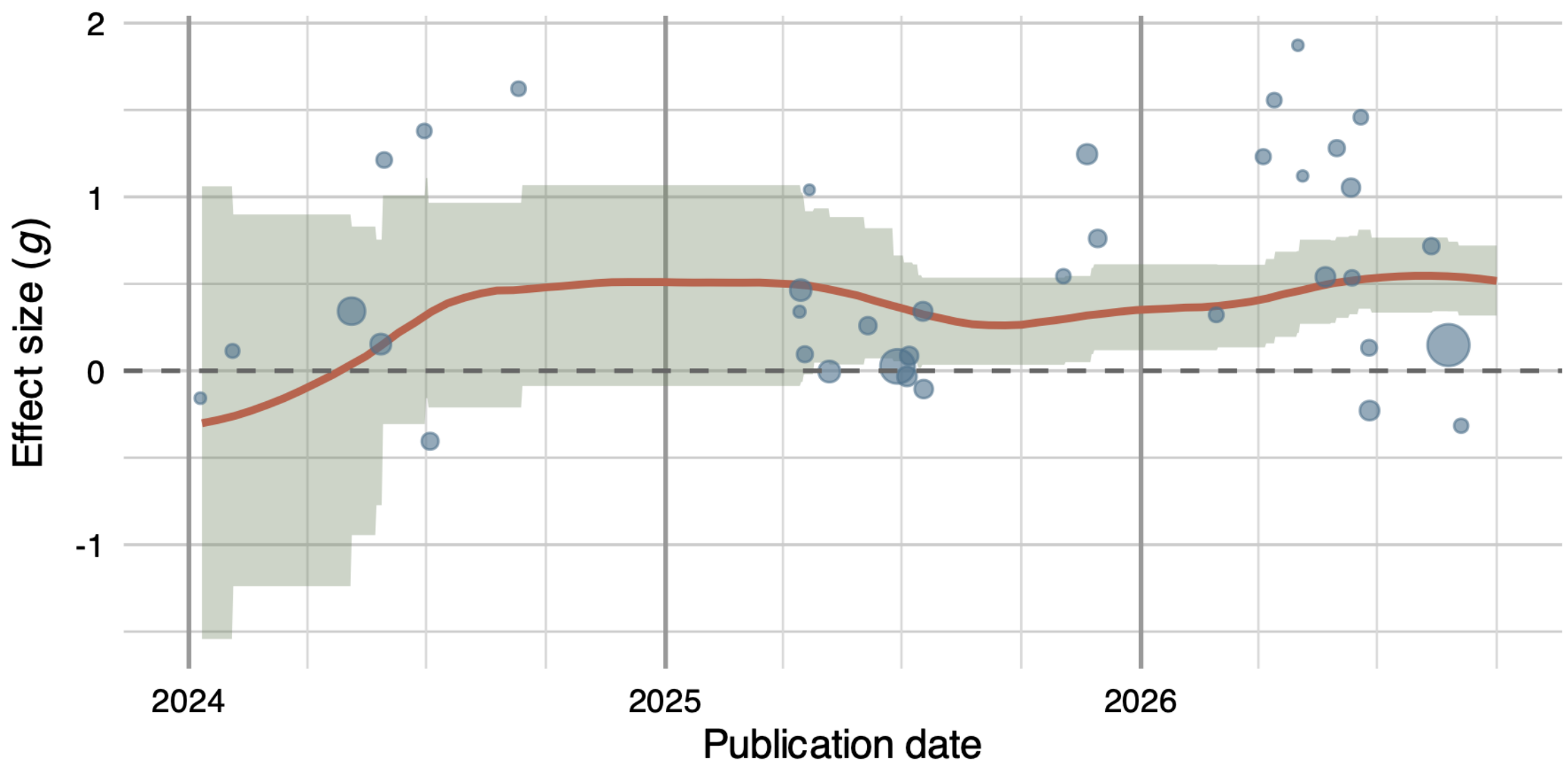

**Figure 2.** Cumulative Bayesian meta-analysis over time.

*Note.* Study-level effect estimates (Hedges' *g*) are shown as points at their publication dates, with point size proportional to the effective sampling precision of each study, accounting for within-study dependence. The smoothed line and shaded region indicate the posterior median and 95% credible interval of the pooled effect as evidence accumulates.

heterogeneity, expressed as the posterior median of the study-level random-effect standard deviation, was estimated at $\tau$ = 0.42 (95% CrI [0.23, 0.64]), suggesting considerable variability across studies with wide uncertainty. Additional variability among multiple effect sizes within studies was also observed ($\tau$ = 0.39, 95% CrI [0.28, 0.54]).

### 4.3 Moderator analyses

The results of the moderator analyses are shown in Table 1. For most moderators, Bayes factors were within the range of *anecdotal evidence* ($0.33 < BF_{10} < 3$), indicating that the amount of evidence is still limited. For *intervention dose* (log-transformed in hours), a Bayes factor of $BF_{10}$ = 67 indicated *strong evidence* for an association with effect size. The estimated effect was $g$ = 0.39 at one hour and increased by $\Delta g$ = 0.14 for each doubling of intervention dose. For *integration mode*, a Bayes factor of $BF_{10}$ = 61 indicated *strong evidence* for including the moderator. Estimated category effects were g = 0.68 when generative AI supplemented regular instruction, $g$ = 0.09 when it replaced teachers, and $g$ = 0.21 when it was used as a standalone tool. Strong evidence also favored including *grouping* (BF10 = 30); category estimates were lower in individual settings ($g$ = 0.41) than in group settings ($g$ = 1.08), teacher-led settings ($g$ = 0.92), or tutoring settings ($g$ = 0.83). Moderate evidence favored including *educational level* and *output mode*; category estimates were higher in primary education and with picture or video output. Prior AI literacy was not assessed by any of the included studies, and prior general abilities were assessed by only three studies. Therefore, neither moderator was included in the analyses. Figure 3 shows how Bayes factors changed across versions.

**Table 1.** Moderator analyses

| | *N* | *k* | *Est.* | 95% CrI | $BF_{10}$ | $\tau^2$(N) | $\tau^2$(k) | $R^2$(N) | $R^2$(k) |
|---|---|---|---|---|---|---|---|---|---|
| **General Intervention Characteristics** | | | | | | | | | |
| *Educational Level*[1] | | | | | 9.889 | 0.090 | 0.170 | .496 | <.001 |
| Primary and lower | 12 | 18 | 0.830 | [0.532, 1.115] | | | | | |
| Secondary | 12 | 39 | 0.328 | [0.063, 0.629] | | | | | |
| Post-Secondary, Tertiary | 11 | 15 | 0.404 | [0.119, 0.698] | | | | | |
| Adults | 3 | 3 | 0.393 | [-0.022, 0.824] | | | | | |
| *Mathematical Content*[1] | | | | | 0.696 | 0.163 | 0.160 | .090 | <.001 |
| Number and Operations | 10 | 16 | 0.597 | [0.302, 0.894] | | | | | |
| Algebra | 24 | 50 | 0.442 | [0.216, 0.678] | | | | | |
| Geometry | 9 | 19 | 0.574 | [0.183, 0.957] | | | | | |
| Measurement | 6 | 13 | 0.730 | [0.270, 1.178] | | | | | |
| Data Analysis and Probability | 5 | 9 | 0.498 | [0.104, 0.894] | | | | | |
| *Grouping*[1] | | | | | 30.006 | 0.103 | 0.149 | .440 | .036 |
| Individual | 32 | 69 | 0.408 | [0.222, 0.613] | | | | | |
| Group | 7 | 9 | 1.078 | [0.585, 1.548] | | | | | |
| Teacher-led plenary | 8 | 11 | 0.918 | [0.451, 1.380] | | | | | |
| Tutoring | 1 | 1 | 0.833 | [0.296, 1.371] | | | | | |
| *Intervention Dose* | | | | | 66.854 | 0.105 | 0.149 | .416 | .034 |
| At 1 hour | 37 | 74 | 0.389 | [0.204, 0.582] | | | | | |
| Doubling | 37 | 74 | 0.141 | [0.061, 0.221] | | | | | |
| **AI-Specific Intervention Characteristics** | | | | | | | | | |
| *Purpose of AI Use*[1] | | | | | 0.545 | 0.185 | 0.155 | <.001 | <.001 |
| Mathematics Expert | 5 | 8 | 0.374 | [-0.059, 0.809] | | | | | |
| Dialogic Partner | 25 | 56 | 0.542 | [0.302, 0.783] | | | | | |
| Material Design | 3 | 6 | 0.520 | [0.091, 0.950] | | | | | |
| Facilitator | - | - | - | - | | | | | |
| Teacher Support | 9 | 11 | 0.530 | [0.154, 0.902] | | | | | |
| *Modified AI* | | | | | 0.356 | 0.177 | 0.159 | .010 | <.001 |
| No | 9 | 23 | 0.540 | [0.252, 0.830] | | | | | |
| Yes | 30 | 51 | 0.513 | [0.311, 0.724] | | | | | |
| *Output Mode*[1] | | | | | 4.175 | 0.134 | 0.158 | .253 | <.001 |
| Text | 37 | 74 | 0.449 | [0.257, 0.657] | | | | | |
| Pictures | 8 | 12 | 0.959 | [0.425, 1.482] | | | | | |
| Video | 5 | 7 | 0.812 | [0.241, 1.385] | | | | | |
| State of Technology | | | | | 1.655 | 0.152 | 0.154 | .152 | <.001 |
| 2021 | 1 | 1 | 0.456 | [-0.021, 0.932] | | | | | |
| 2022 | 6 | 11 | 0.420 | [0.075, 0.775] | | | | | |
| 2023 | 10 | 33 | 0.524 | [0.225, 0.838] | | | | | |
| 2024 | 11 | 16 | 0.394 | [0.107, 0.699] | | | | | |
| 2025 | 8 | 12 | 0.741 | [0.404, 1.064] | | | | | |
| 2026 | 1 | 1 | 0.613 | [0.134, 1.087] | | | | | |
| *Integration Mode* | | | | | 60.636 | 0.161 | 0.128 | .099 | .167 |
| Replace | 4 | 5 | 0.094 | [-0.351, 0.521] | | | | | |
| Standalone | 8 | 32 | 0.206 | [-0.112, 0.545] | | | | | |
| Supplement | 27 | 37 | 0.676 | [0.463, 0.894] | | | | | |
| Supported | - | - | - | - | | | | | |
| **Extrinsic Characteristics** | | | | | | | | | |
| *Publication Type* | | | | | 0.406 | 0.191 | 0.153 | <.001 | .004 |
| Quartile 1 journal | 16 | 46 | 0.567 | [0.305, 0.837] | | | | | |
| Journal | 6 | 7 | 0.524 | [0.147, 0.903] | | | | | |
| Book Chapter | 1 | 1 | 0.574 | [0.067, 1.078] | | | | | |
| Proceedings | 11 | 14 | 0.453 | [0.148, 0.764] | | | | | |
| Preprint | 3 | 6 | 0.493 | [0.065, 0.923] | | | | | |
| **Methodological Characteristics** | | | | | | | | | |
| *Research Design* | | | | | 1.392 | 0.172 | 0.147 | .038 | .046 |
| Randomized controlled trial | 12 | 14 | 0.350 | [0.050, 0.665] | | | | | |
| Cluster-randomized trial | 8 | 35 | 0.702 | [0.365, 1.042] | | | | | |
| Non-randomized study | 17 | 25 | 0.540 | [0.275, 0.812] | | | | | |
| *Assessment Timepoint* | | | | | 0.882 | 0.170 | 0.159 | .048 | <.001 |
| Immediate | 37 | 71 | 0.509 | [0.312, 0.712] | | | | | |
| Delayed | 3 | 3 | 0.715 | [0.217, 1.218] | | | | | |
| *Research Setting* | | | | | 0.545 | 0.185 | 0.156 | <.001 | <.001 |
| Field | 26 | 59 | 0.511 | [0.286, 0.741] | | | | | |
| Controlled | 10 | 14 | 0.538 | [0.194, 0.877] | | | | | |
| Lab | 1 | 1 | 0.601 | [0.021, 1.176] | | | | | |

*Note.* [1]For these moderators, several levels could be observed within the same effect. In these cases, the effect was distributed evenly across the levels by incremental dummy coding. Thus, counts can sum to more than the totals. N = studies; k = effects; Est. = estimated Hedges' g (change in g per doubling for dose); $\tau^2$ = residual heterogeneity variance; $R^2$ = proportional variance reduction.

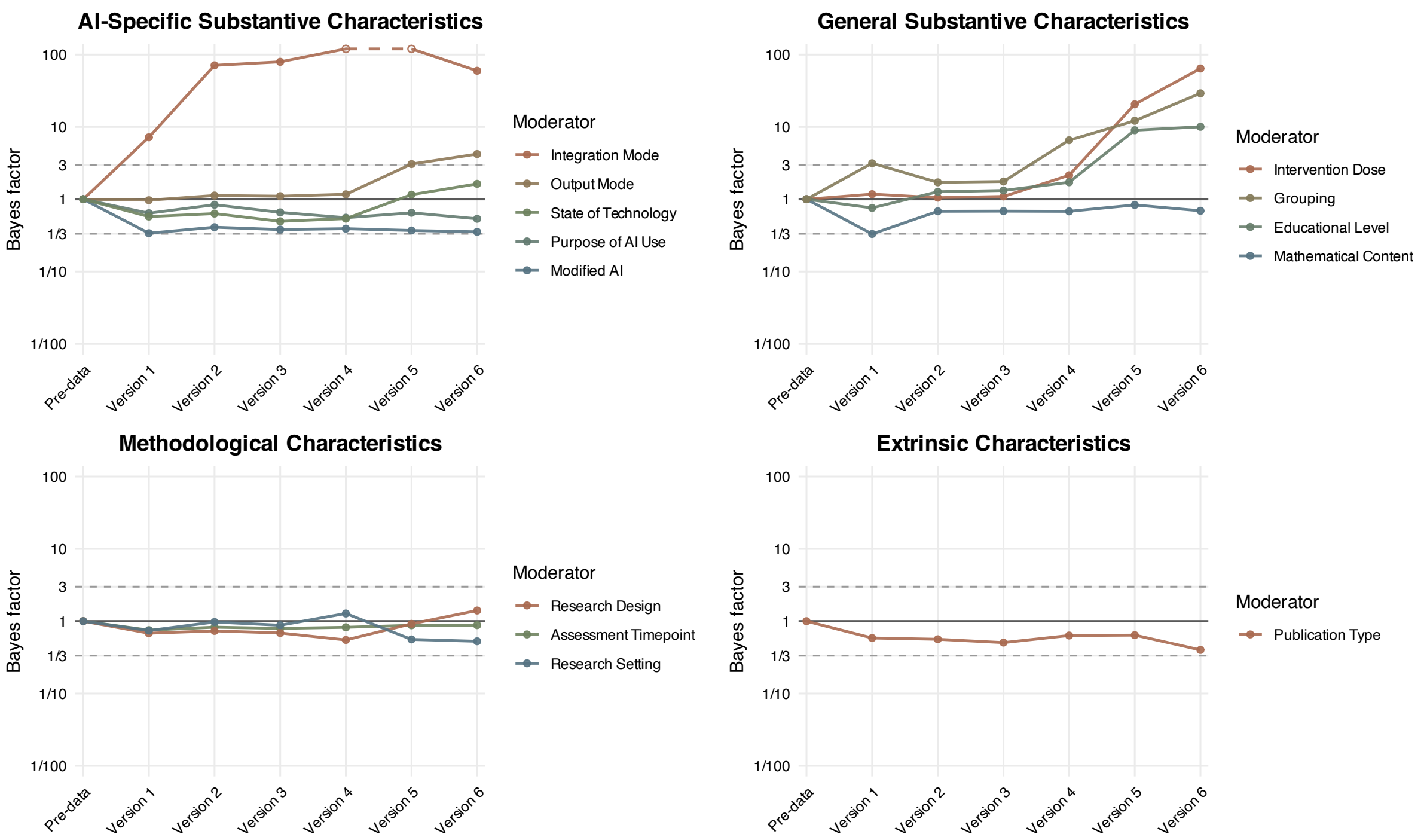

**Figure 3.** Bayes factors for moderator effects across successive versions, grouped by moderator category.

*Note.* Panels distinguish general intervention, AI-specific intervention, methodological, and extrinsic characteristics. Colored lines show individual moderators. Bayes factors above 1 favor moderator inclusion, whereas values below 1 favor exclusion. The pre-data value of 1 represents equal support for inclusion and exclusion. Dashed horizontal lines mark BF = 3 and BF = 1/3, conventional thresholds for moderate evidence for inclusion and exclusion, respectively (Lee & Wagenmakers, 2014). The y-axis is logarithmic. For display, values are capped at 1/120 and 120. Open circles indicate values above 120; dashed line segments connect successive values exceeding this upper limit.

## 4.4 Publication Bias

Publication bias was assessed using the multilevel robust Bayesian model-averaged framework implemented in *RoBMA* (Bartoš & Maier, 2020; Bartoš, Maier, et al., 2025). The multilevel specification accounted for clustering of multiple effects within studies, but not for correlated sampling errors because *RoBMA* does not currently allow a full sampling variance–covariance matrix. The standard model ensemble combined six selection-model weight functions with PET–PEESE. This analysis indicated anecdotal evidence against publication bias (inclusion $BF$ = 0.36; posterior probability = .27). The corresponding model-averaged effect estimate was $g$ = 0.50, 95% CrI [0.00, 0.72], with strong evidence for the presence of an overall effect ($BF_{10}$ = 26.06; posterior probability = .96). The model-averaged funnel plot is given in Figure 4.

## 4.5 Changes to the results

Compared with Version 5, we added five studies contributing nine effects. As in the previous version, edge cases in the newly screened studies required us to make our inclusion criteria more explicit and stricter regarding AI use in control groups. This led to the exclusion of one study that had been included in Version 5 (Luo et al., 2026). One study (Dubey et al., 2026) was excluded after we were made aware that it reported data from a paper by another group of researchers without

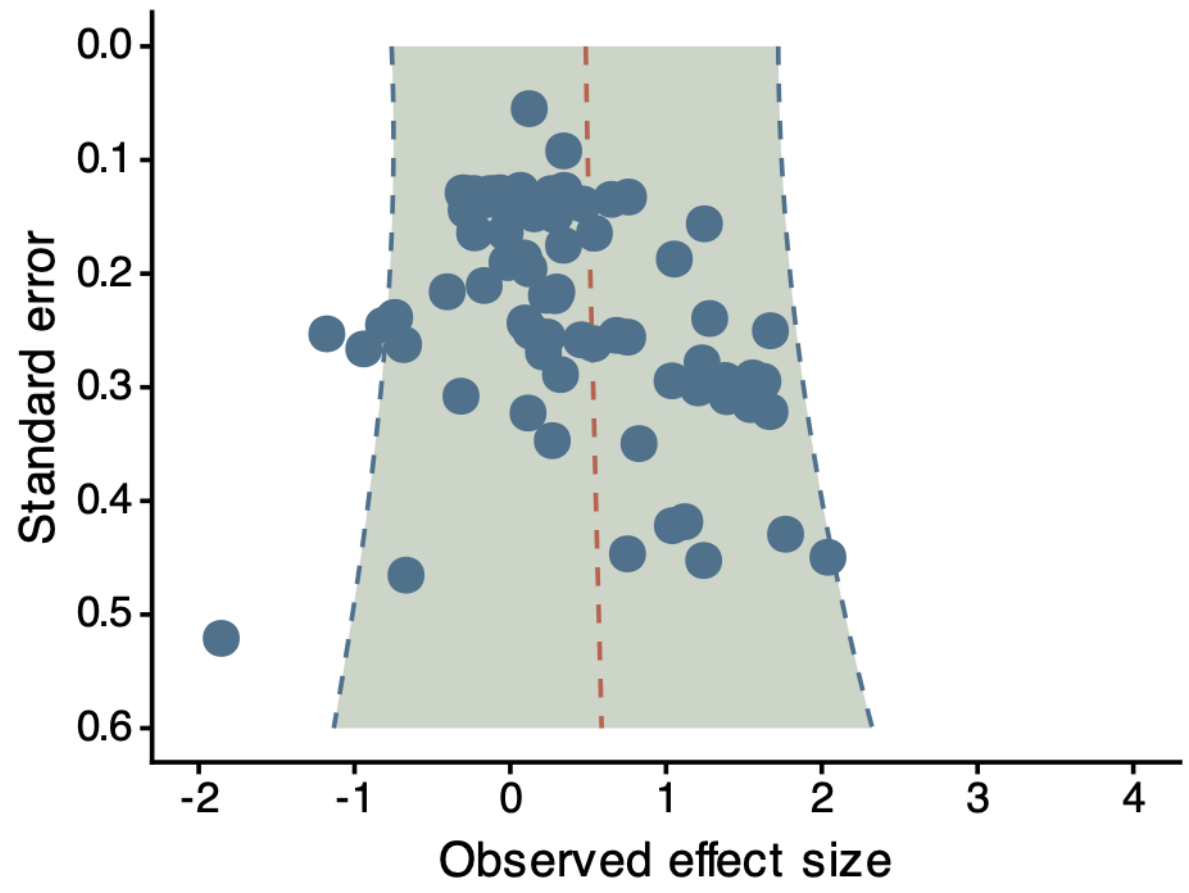

**Figure 4.** Model-averaged funnel plot for the publication-bias analysis.

*Note.* Observed effect sizes are plotted against their standard errors. The funnel boundaries represent the central 95% region of the sampling distribution implied by the fitted *RoBMA* model, incorporating heterogeneity and publication-bias adjustment.

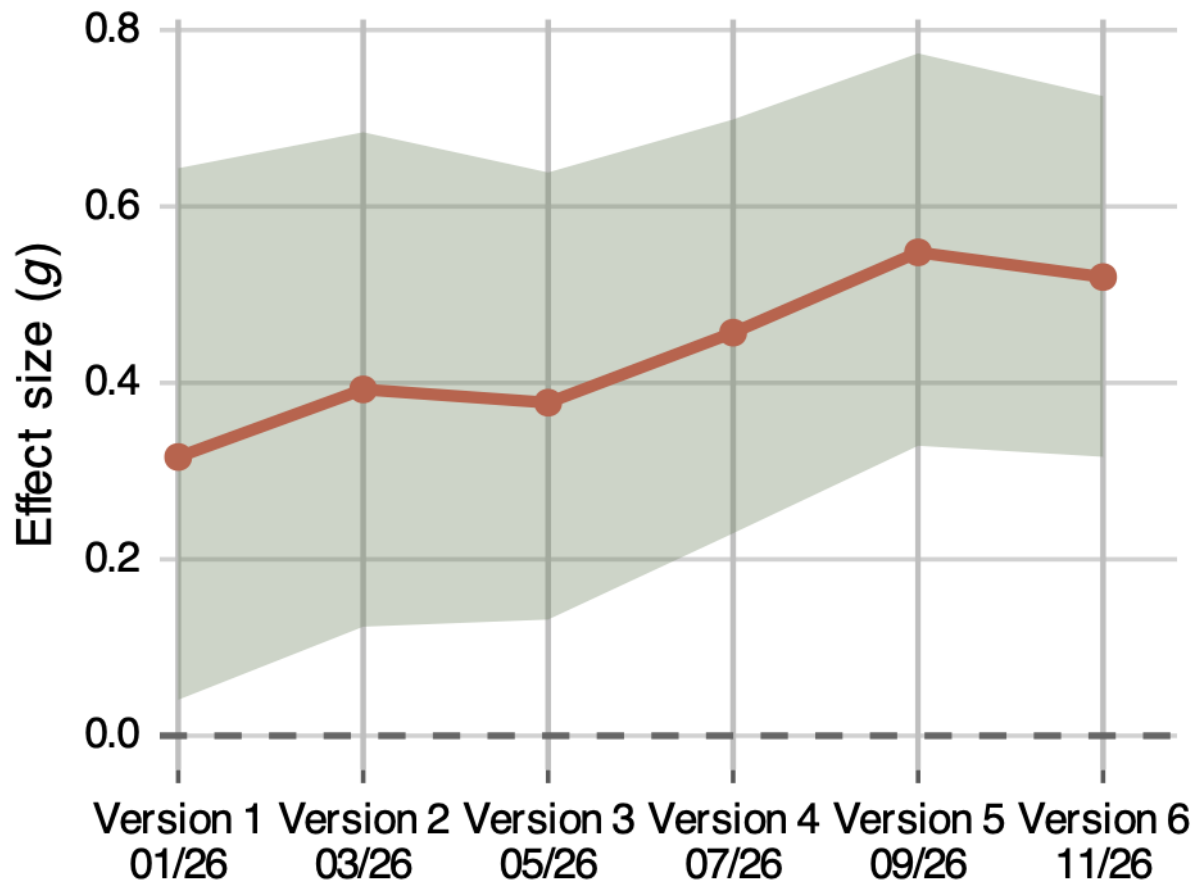

**Figure 5.** Cumulative pooled effect estimates across versions.

*Note.* Points indicate the posterior median Hedges' *g* for each version, and shaded bands represent the corresponding 95% credible intervals.

proper acknowledgement. After correspondence with the corresponding authors of both publications, we decided to exclude Dubey et al. (2026) from future versions.

Overall, the number of included effects therefore increased by five. The posterior median pooled effect decreased by $\Delta g$ = 0.05. We included four new moderators and report *intervention dose* instead of *intervention span*. *Grouping* and *intervention dose* showed strong evidence in favor of inclusion. The Bayes factor for *integration mode* decreased but is still strongly in favor of the moderator. The Bayes factor for *educational level* increased from 5.11 to 9.89, remaining within the range of *moderate evidence.* The publication-bias analysis again provided anecdotal evidence against bias. Prior sensitivity results are reported in Section 3.5.

## 5 Discussion

LLAMA LIMA provides an ongoing synthesis of intervention studies that use generative AI to support mathematics learning. Our analysis shows a positive median effect ($g$ = 0.52, 95% CrI [0.32, 0.72]) across 37 studies, providing evidence for an overall positive effect. Since the first version was published in January 2026, this effect has increased (Figure 5). Whether this reflects an increasing potential of generative AI to support mathematics learning remains unclear, since evidence for an association with state of technology remained inconclusive.

To situate the size of the effect, meta-analytic estimates for other interventions in mathematics education (e.g., digital media $g$ = 0.55; Hillmayr et al., 2020; visualizations $g$ = 0.50; Schoenherr et al., 2024) can be used, suggesting that the current pooled effect is similar in magnitude to estimates for other systematic interventions. Specifically for generative AI, Bartoš, Martinková, et al. (2025) report a smaller effect of $g$ = 0.37 for ChatGPT use across educational domains when correcting for publication bias. Similarly, Boolzen et al. (2026) report a substantially smaller effect across STEM domains when taking publication bias into account ($g$ = 0.076). Some other meta-analyses on generative AI report substantially higher effect sizes, but considerable reservations about their validity and trustworthiness have been documented in detail by O'Neill (2026) and Fütterer et al. (2026).

One of the issues not properly accounted for in many existing meta-analyses in the field is publication bias (O'Neill, 2026). Analyses by Bartoš, Martinková, et al. (2025) and Boolzen et al. (2026) indicate substantial publication bias in broader educational and STEM research. In contrast, our analyses, using the RoBMA framework (Bartoš & Maier, 2020), provided anecdotal evidence against publication bias in the current dataset. This difference may be related to our specific focus on mathematics education and our study eligibility criteria, although further research is needed to establish whether publication bias indeed differs across educational domains.

The substantial heterogeneity of effects across studies suggests that the effectiveness of generative AI in mathematics education varies across implementations and contexts. The current moderator analyses included 13 moderators. We found strong evidence in favor of including three moderators and moderate evidence for including two; evidence for the remaining moderators was anecdotal.

Strong evidence favored including integration mode. The estimated effects were larger when generative AI supplemented regular instruction than when it replaced teachers ($\Delta g$ = 0.58). This finding supports the view that generative AI should be integrated in ways that combine the strengths of teachers and generative AI (Giannakos et al., 2025) and is in line with findings by Boolzen et al. (2026). Teacher guidance may help align AI use with the learning goal, elicit mathematical explanations, and check outputs, while supporting learners' own cognitive activity (Cho, 2026; Reinhold et al., 2024). Cho's (2026) framework offers a possible theoretical interpretation: teacher mediation and task design may help preserve learners' own generative activity. Our moderator analyses do not directly test this mechanism.

Strong evidence was also found for an association with intervention dose, such that greater total session time was associated with larger effects, and for grouping, with smaller effects in individual settings. Moderate evidence favored including educational level and output mode; estimated category effects were higher in primary education and with picture or video output. Evidence for an association with state of technology remained inconclusive. For the purpose of AI use, the Bayes factor provided only anecdotal evidence against an influence of the moderator. The same purpose of AI use may involve different learner activities, suggesting that purpose categories alone may not capture the cognitive processes underlying learning (Cho, 2026). For research setting, AI modification, and mathematical content area, there was anecdotal evidence against moderator inclusion. None of the studies included in our analyses assessed students' AI-literacy, and only three studies assessed general academic abilities. Given the importance of students' prior abilities for learning, this indicates a substantial blind spot in the current research landscape.

With the increasing number of studies, recurring issues of trustworthiness, and the relatively short delay between publication and inclusion in our analysis, LLAMA LIMA arguably demonstrates the value of a living meta-analytic approach for synthesizing evidence in rapidly evolving research fields. Living meta-analyses have been used in other disciplines but remain rare in education. Howard and Slemp (2026) developed a living meta-analytic database without publishing successive versions of its analyses. To our knowledge, LLAMA LIMA is the first living meta-analysis in education to publish regularly updated versions of the synthesis, following the approach to living evidence reporting described by Elliott et al. (2014) and Akl et al. (2024).

The changes made in Version 6 also illustrate the importance of correction and accountability. One study was excluded because we had incorrectly coded its control group, and another was removed after correspondence confirmed that it used data from a different publication without acknowledging this in the paper (Dubey et al., 2026). These cases underline the concerns about traceability and trustworthiness raised by Fütterer et al. (2026) and O'Neill (2026). The living format makes such decisions and their consequences transparent and keeps the synthesis open to feedback, critique, and correction.

### 5.1 Limitations

The results of this versioned living meta-analysis are provisional and may change as new studies become available. In addition, capabilities, interfaces, and usability of generative AI evolve rapidly, such that effects might not be comparable across longer periods of time. We address this limitation by making updates and versions transparent, following recommendations for this publication format (Akl et al., 2024). The inclusion of preprints reduces publication lag but potentially includes evidence that has not been subject to rigorous peer review. However, in the current version, three included studies are available only as preprints, and other studies that had previously been published as preprints have now been retrieved from journals or published conference proceedings. Conclusions are limited by the quality of the primary studies and variation in interventions, outcome measures, and assessment timepoints. Moderator associations may reflect confounding between study characteristics and do not establish causal effects.

### 5.2 Conclusion

Our analyses indicate that generative AI systems may support mathematics learning, particularly when integrated with regular teaching and during group work or teacher-led or tutor-led activities. At the same time, the currently available empirical evidence on effective implementations, their prerequisites, and contextual conditions remains limited relative to the breadth of expectations and proposed applications discussed in the literature and public discourse. Accordingly, further rigorous research is needed, alongside ongoing and systematic evaluation, to clarify whether and under which conditions generative AI can reliably contribute to mathematics learning. With LLAMA LIMA, we continuously observe and evaluate this development.

## Appendix A.1

Boolean search string (Scopus):

("generative AI" OR "generative Artificial Intelligence" OR genAI* OR "large language model*" OR LLM* OR "chatbot*" OR "AI tutor*" OR "AI assistant*" OR ChatGPT OR "Chat GPT") AND (educat* OR teach* OR instruct* OR pedagogy OR curriculum OR classroom OR student* OR school* OR "higher education" OR "K-12") AND (math* OR algebra OR geometry OR arithmetic OR calculus OR "word problem*" OR numeracy OR fraction* OR "quantitative reasoning") AND (intervention* OR experiment* OR "control group" OR RCT OR "randomized control*" OR treatment OR pretest OR posttest OR (pre* AND post*))

Semantic search string (OpenAlex):

"Experimental or quasi-experimental intervention study on generative artificial intelligence in mathematics education, where students or teachers use ChatGPT, a large language model, generative AI chatbot, or AI tutor, and mathematics learning outcomes are compared with a control group without generative AI."

**Appendix Table A.2.** Overview of the included studies with selected study features.

| Authors | Publication date | State of technology | Publication format | Number of participants | Educational level | Content area | AI purpose | Integration mode | AI system modification | Number of effects | Pooled study effect (*g*) |
|---|---|---|---|---|---|---|---|---|---|---|---|
| Aslan et al. (2026) | 17.06.2026 | 2025 | JQ1 | 51 | primary | A | 2;5 | spl | yes | 2 | 1.46 |
| Bastani et al. (2025) | 25.06.2025 | 2023 | JQ1 | 943 | secondary | N; A; G; D | 2 | sta | yes; no | 24 | 0.03 |
| Bešlić et al. (2024) | 03.07.2024 | 2024 | proc. | 86 | secondary | N | 2 | spl | yes | 1 | -0.40 |
| Canonigo (2024) | 13.09.2024 | 2023 | JQ1 | 60 | secondary | A | 2 | spl | no | 1 | 1.62 |
| Chen, Judicke, et al. (2026) | 27.06.2026 | 2024 | proc. | 63 | post-secondary | A | 1;2 | spl | yes | 1 | -0.23 |
| Chen, Tang, et al. (2026) | 27.06.2026 | 2025 | proc. | 148 | adults | A | 2 | sta | yes | 1 | 0.13 |
| Dasari et al. (2024) | 11.01.2024 | 2022 | JQ1 | 29 | post-secondary | A | 1 | rep; spl | no | 2 | -0.16 |
| El-Shara et al. (2025) | 18.04.2025 | 2023 | jour. | 66 | post-secondary | A | 2 | spl | yes | 1 | 0.09 |
| **Fekumo and Phoshoko (2026)** | **27.05.2026** | **2025** | **jour.** | **163** | **post-secondary** | **A** | **1;2** | **spl** | **no** | **1** | **0.54** |
| Fardian et al. (2025) | 15.04.2025 | 2022 | JQ1 | 30 | post-secondary | A | 1 | rep; spl | no | 2 | 0.34 |
| Güney et al. (2027) | 03.05.2026 | 2025 | JQ1 | 29 | primary | M | 2;5 | spl | yes | 2 | 1.87 |
| Hau et al. (2026) | 07.04.2026 | 2023 | JQ1 | 60 | primary | G;M | 5 | spl | yes | 1 | 1.23 |
| Henkel et al. (2024) | 05.05.2024 | 2023 | proc. | 477 | primary; secondary | N;A | 2 | sta | yes | 1 | 0.34 |
| Karaman (2024) | 04.02.2024 | 2023 | JQ1 | 39 | primary | G;M | 5 | spl | yes | 1 | 0.11 |
| Khaled and Alghfeli (2026) | 14.04.2026 | 2025 | pre. | 60 | secondary | A | 2 | spl | no | 1 | 1.56 |
| Koronaiou and Kostas (2026) | 01.05.2026 | 2024 | chap. | 25 | primary | N | 5 | spl | yes | 1 | 1.12 |
| Lademann et al. (2025) | 07.05.2025 | 2023 | JQ1 | 214 | secondary | N | 5 | sta | yes | 1 | 0.00 |
| Li and Lyu (2025) | 03.07.2025 | 2021 | JQ1 | 151 | post-secondary | A | 2 | sta | yes | 1 | -0.03 |
| Liu et al. (2025) | 25.04.2025 | 2024 | proc. | 24 | primary | G | 2 | spl | yes | 1 | 1.04 |
| Malik et al. (2026) | 01.03.2026 | 2025 | jour. | 47 | secondary | A | 2 | spl | yes | 1 | 0.32 |
| **Miao et al. (2026)** | **25.08.2026** | **2024** | **pre.** | **1628** | **secondary** | **N; A; G; M; D** | **2** | **spl** | **yes** | **4** | **0.15** |
| Nieto-Jalil et al. (2026) | 10.06.2026 | 2024 | proc. | 135 | post-secondary | A | 2 | spl | yes | 1 | 1.05 |
| Pardos and Bhandari (2024) | 24.05.2024 | 2022 | jour. | 184 | adults | A;D | 5 | sta | yes | 1 | 0.15 |
| **Rief et al. (2026)** | **02.09.2026** | **2024** | **pre.** | **42** | **secondary** | **A** | **2;5** | **spl** | **yes** | **1** | **-0.32** |
| Rücker and Becker-Genschow (2025) | 24.11.2025 | 2023 | JQ1 | 195 | secondary | N | 2 | spl | yes | 1 | 1.25 |
| Serrano Heredia et al. (2025) | 03.06.2025 | 2024 | proc. | 88 | post-secondary | A;G;M | 3 | rep | yes | 2 | 0.26 |
| Steinbach et al. (2025) | 17.07.2025 | 2022 | proc. | 131 | adults | A;D | 5 | sta | yes | 1 | 0.34 |
| Steindl et al. (2025) | 04.11.2025 | 2024 | proc. | 49 | post-secondary | A | 1 | sta | no; yes | 2 | 0.54 |
| Sánchez et al. (2026) | 10.06.2026 | 2024 | proc. | 59 | post-secondary | A | 2 | spl | yes | 1 | 0.53 |
| **Tasdelen and Bodemer (2025)** | **10.07.2025** | **2024** | **JQ1** | **114** | **primary** | **N** | **3** | **spl** | **yes** | **1** | **0.09** |
| Utami et al. (2024) | 29.05.2024 | 2022 | JQ1 | 51 | primary | G;M | 3 | spl | yes | 3 | 1.21 |
| Wahba et al. (2024) | 01.07.2024 | 2023 | jour. | 56 | post-secondary | D | 2 | spl | no | 1 | 1.38 |
| Wang and Kuo (2026) | 31.05.2026 | 2026 | JQ1 | 83 | primary | N | 2 | rep | yes | 1 | 1.28 |
| Xie et al. (2025) | 26.11.2025 | 2025 | proc. | 85 | primary | G | 2 | spl | no | 2 | 0.76 |
| Xing et al. (2025) | 17.04.2025 | 2023 | JQ1 | 212 | secondary | A;G | 2 | spl | yes | 1 | 0.46 |
| Xuan et al. (2025) | 14.07.2025 | 2022 | JQ1 | 120 | secondary | A | 2 | spl | yes | 2 | -0.11 |
| **Yang et al. (2026)** | **13.08.2026** | **2025** | **jour.** | **64** | **primary** | **N; A** | **2** | **spl** | **yes** | **2** | **0.72** |

*Note*. Newly included studies are in **bold**. Publication format: chap. = Book chapter; jour. = Journal; JQ1 = Journal ranked in the 1st quartile according to WoS; pre. = Preprint; proc. = Proceedings; Content areas (NCTM): N = Number and Operations, A = Algebra, G = Geometry, M = Measurement, D = Data Analysis and Probability; AI purpose: 1 = Direct mathematical support for learners (Mathematics Expert); 2 = Learner–AI dialogue (Dialogic Partner); 3 = Generation of learning materials (Material Design); 4 = Facilitation of collaborative learning (Facilitator); 5 = Teacher- or researcher-mediated support (Teacher Support). Integration Mode: sta = standalone; spl = supplement; rep = replacement, spr = supported.